\documentclass[12pt]{article}
\usepackage{amsmath,amssymb}
\usepackage[final]{epsfig}

\newtheorem{theorem}{Theorem}[section]

\newtheorem{theorem-definition}[theorem]{Theorem-Definition}
\newtheorem{theorem-construction}[theorem]{Theorem-Construction}
\newtheorem{lemma-definition}[theorem]{Lemma--Definition}
\newtheorem{lemma}[theorem]{Lemma}
\newtheorem{proposition}[theorem]{Proposition}

\newtheorem{conjecture}[theorem]{Conjecture}
\newtheorem{definition}[theorem]{Definition}

\begin{document}
\newcommand{\Z}{{\Bbb Z}}
\newcommand{\R}{{\Bbb R}}
\newcommand{\Q}{{\Bbb Q}}
\newcommand{\C}{{\Bbb C}}
\newcommand{\lms}{\longmapsto}
\newcommand{\lra}{\longrightarrow}
\newcommand{\hra}{\hookrightarrow}
\newcommand{\ra}{\rightarrow}

\begin{titlepage}
\title{Euclidean scissor congruence groups and mixed 
Tate motives over dual numbers}  
\author{A. B. Goncharov}
\end{titlepage}
\date{}
\maketitle

     \qquad \qquad \qquad {\it To Spencer Bloch, with admiration,  for his 60th birthday}

\tableofcontents

 \section{Introduction} 

{\bf Summary}. We define Euclidean scissor congruence groups 
for an arbitrary algebraically closed field $F$ and  formulate a  
conjecture describing 
them. Using the Euclidean and Non-Euclidean $F$--scissor congruence groups 
we construct a category which is  conjecturally equivalent to a 
subcategory of the category ${\cal M}_T(F_{\varepsilon})$ of mixed Tate motives 
over the dual numbers $F_{\varepsilon}:= F[\varepsilon]/\varepsilon^2$.  

{\bf 1. Euclidean scissor congruence groups and 
a generalization of  Hilbert's third problem}. 
Let $F$ be an arbitrary algebraically closed field. 
In Chapter 3 of 
\cite{Gvol} 
we defined an $F$--scissor congruence 
group $S_n(F)$ of polyhedrons in the projective space $P^{2n-1}(F)$  
equipped with a non-degenerate quadric $Q$. 
The classical spherical and hyperbolic scissor congruence 
groups are subgroups of $S_n(\C)$. The direct sum 
$$
S_{\bullet}(F):= \oplus_{n \geq 0}S_{n}(F); \qquad S_0(F) =\Q
$$
is equipped with a structure of a commutative, graded Hopf $\Q$--algebra. 
The coproduct is given by the Dehn invariant map
$$
D: S_{n}(F) \lra \oplus_{0 \leq k\leq n} S_{k}(F)\otimes S_{n-k}(F)
$$

In this paper we define 
{\it Euclidean $F$--scissor congruence groups} ${\cal E}_n(F)$ of polyhedrons 
in $(2n-1)$--dimensional affine space over $F$ equipped 
with a non-degenerate quadratic form $Q$. If $F=\R$ (which is  not algebraically closed!) 
and $Q$ is positive definite, we get the classical Euclidean scissor 
congruence group ${\cal E}_n(\R)$ in $\R^{2n-1}$. We define the 
Dehn invariant map
$$
D^E: {\cal E}_{n}(F) \lra \oplus_{1 \leq k\leq n} {\cal E}_{k}(F)\otimes S_{n-k}(F)
$$
and show that it provides the graded $\Q$--vector space 
$$
{\cal E}_{\bullet}(F):= \oplus_{k \geq 0}{\cal E}_{n}(F)
$$
with a structure of a comodule over the Hopf algebra $S_{\bullet}(F)$. 
The cobar complex calculating the cohomology of this comodule 
looks as follows:
\begin{equation} \label{1.15.03.12xy}
{\cal E}_{\bullet}(F) \stackrel{}{\lra}  
{\cal E}_{\bullet}(F)\otimes {S}_{\bullet}(F) 
\stackrel{}{\lra} {\cal E}_{\bullet}(F)\otimes {S}_{\bullet}(F)^{ \otimes 2} \lra ... 
\end{equation}
The differential is cooked up from $D^E$ and $D$ using the Leibniz rule. 
For example the first two arrows are $D^E$ and 
$D^E \otimes {\rm Id} - {\rm Id} \otimes D$. 

We place the first group in degree $1$, 
denote the complex by  ${\cal E}^*_{(\bullet)}(F)$,
 and call it {\it Euclidean Dehn complex}. 
The differential preserves the grading. We denote by ${\cal E}^*_{(n)}(F)$ 
the degree $n$ subcomplex of (\ref{1.15.03.12xy}). 

The action by dilotations of the group $F^*$ in a Euclidean  $F$-vector space   
provides an $F^*$--action  on  the Euclidean $F$--scissor congruence groups. 
So (\ref{1.15.03.12xy}) a complex of $F^*$--modules. 
For  an $F$--vector space $V$ denote by $V<p>$ 
the twisted $F^*$--module structure $*$ on $V$ 
given by $f * v = f^{2p+1}\cdot v$. 
The volume provides a homomorphism 
$$
{\rm Vol}: {\cal E}_{n}(F) \lra F<n>
$$
\begin{conjecture} \label{1.15.03.12} 
a) Let $F$ be an arbitrary algebraically closed field. Then we have 
canonical isomorphism 
of $F^*$--modules
$$
H^i({\cal E}_{(n)}^{\ast}(F)) = \Omega^{i-1}_{F/\Q}<n>
$$

b) The same is true in the classical case $F=\R$. 
\end{conjecture}
One can view this as a generalization of 
Hilbert's Third Problem. Indeed, according to Sydler's theorem \cite{S} 
the following complex is exact:
$$
0 \lra \R \lra {\cal E}_2(\R) \stackrel{D^E}{\lra} \R \otimes S^1 \lra \Omega^1_{\R/\Q}
\lra 0
$$
In particular 
the kernel of the 
Dehn invariant is identified by the volume 
homomorphism with 
$\R$. 
So the Dehn invariant and the volume 
determine a polyhedron in $\R^3$ uniquely up to scissor congruence.  

Sydler's theorem   gives the $n=2$ case of the part b) of conjecture. 
The $n=2$ case of the part a) can probably be deduced 
from the results of 
Sydler, Dupont, Cathelineau, Sah in \cite{S}, \cite{D}, \cite{DS}, \cite{Ca}.

Another  key problem  is 
the structure of the groups ${\cal E}_{n}(F)$. 
By the Euler characteristic argument (and induction) the answer is controlled, 
although in a cryptic way,  by 
Conjecture \ref{1.15.03.12}. Let us get it in a more explicit form.

{\bf 2. A hypothetical description of the Euclidean scissor congruence groups.}
Let 
$$
{\cal Q}_{\bullet}(F):= \frac{S_{\bullet}(F)}
{S_{>0}(F) \cdot S_{>0}(F)}
$$
be the space of indecomposables of the Hopf algebra $S_{\bullet}(F)$.
It is a graded Lie coalgebra with the cobracket induced by the coproduct in 
 $S_{\bullet}(F)$. 

Let  
$\Q_{\varepsilon}$--mod be  a tensor category 
with the objects $V = V_0 \oplus V_1 \varepsilon$, where 
$V_0, V_1$ are $\Q$--vector spaces, and ${\rm Hom}$ and tensor product defined by 
$$
{\rm Hom}_{\Q_{\varepsilon}-mod}(V_0 \oplus V_1 \varepsilon, 
W_0 \oplus W_1 \varepsilon) := {\rm Hom}(V_0, W_0) \oplus  {\rm Hom}(V_1, W_1)
$$
$$
(V_0 \oplus V_1 \varepsilon) \otimes_{\Q_{\varepsilon}-mod} (W_0 \oplus W_1 \varepsilon): 
= V_0\otimes W_0 \oplus \Bigl(V_0\otimes W_1 \oplus W_0\otimes V_1\Bigr) \varepsilon
$$
Since $V_1\varepsilon \otimes V_2\varepsilon =0$, it is not a rigid tensor category. 

Observe that a Lie coalgebra ${\cal L}_{\varepsilon}$  in the category 
$\Q_{\varepsilon}$--mod is just the same thing as 
a Lie coalgebra ${\cal L}$ and a comodule ${\cal L}^a$ over it:   
${\cal L}_{\varepsilon} = {\cal L}_{}\oplus {\cal L}^a\cdot \varepsilon$. 

Recall that the Euclidean $F$--scissor congruence groups 
are organized into a comodule ${\cal E}_{\bullet}(F)$ over $S_{\bullet}(F)$, 
and hence over 
${\cal Q}_{\bullet}(F)$. Therefore combining the Euclidean and Non-Euclidean $F$--scissor 
congruence groups we get a  Lie coalgebra 
$$
{\cal Q}_{\bullet}(F_{\varepsilon}):= 
{\cal Q}_{\bullet}(F) \oplus {\cal E}_{\bullet}(F)\cdot \varepsilon
$$
in the category $\Q_{\varepsilon}$--mod. One has ${\cal Q}_{1}(F) =  S_1(F) = F^*$ and ${\cal E}_1(F) = F$.

Recall the higher Bloch groups ${\cal B}_n(F)$ 
 (\cite{G1}-\cite{G2}). We also need their additive versions,  
the $F^*$--modules $\beta_n(F)$, defined 
in  chapter 3 
as extensions of   Cathelineau's groups \cite{Ca}. The $F^*$--module 
$\beta_2(F)$ is 
isomorphic to the one defined  by S. Bloch and H. Esnault  
\cite{BE} in a different way.

Let ${Q}_{\bullet}(F_{\varepsilon})$ be  a 
negatively graded pro-Lie algebra dual to 
${\cal Q}_{\bullet}(F_{\varepsilon})$. 
Denote by ${I}_{\bullet}(F_{\varepsilon})$ its ideal of elements of degree 
$\leq -2$, and by ${\cal I}_{\bullet}(F_{\varepsilon})$ the 
corresponding Lie coalgebra.
Denote by ${\Bbb H}^{\ast}({I}_{\bullet}(F_{\varepsilon}))$ 
the cohomology of ${I}_{\bullet}(F_{\varepsilon})$ in the 
category $\Q_{\varepsilon}$--mod, 
i.e. the cohomology of 
the  standard cochain complex $\Lambda^*{\cal I}_{\bullet}
(F_{\varepsilon})$ in 
the category $\Q_{\varepsilon}$--mod, where 
we take the exterior powers in  $\Q_{\varepsilon}$--mod. 

\begin{conjecture} \label{1.12.03.4qq}  
${I}_{\bullet}(F_{\varepsilon})$ 
 is a free Lie algebra in the category 
$\Q_{\varepsilon}$--mod, with the space of degree $-n$ generators, 
$n=2, 3, ...$, given by 
${\cal B}_{n}(F) \oplus {\beta}_{n}(F)\cdot \varepsilon$. This means that 
$$
{\Bbb H}^p({I}_{\bullet}(F_{\varepsilon})) = 0, \quad p>1; 
\qquad 
{\Bbb H}_{(n)}^1({I}_{\bullet}(F_{\varepsilon})) = 
{\cal B}_{n}(F) \oplus {\beta}_{n}(F)\cdot \varepsilon
$$
\end{conjecture}
Here ${\Bbb H}_{(n)}^1$ stays for the degree $n$ part of ${\Bbb H}^1$. 
In Section 4.2 we will add one more statement  to this conjecture, 
omitted now for the 
sake of simplicity. We would like to stress a similarity 
between this conjecture 
and the freeness conjecture for the mixed elliptic motives, see 
Conjecture 4.3 in \cite{G3}:  
both conjectures can hardly be formulated without use of certain 
esoteric non-rigid tensor structures. 

A statement about an object $V \oplus V_1 \cdot \varepsilon$ in the category  
$\Q_{\varepsilon}$--mod is actually a pair of statements: 
one about $V_0$, and the other about $V_1 \cdot \varepsilon$, called the 
$\Q$- and $\varepsilon$-parts of the statement. 

The $\varepsilon$--part of Conjecture \ref{1.12.03.4qq} is a 
sophisticated version of  
Conjecture \ref{1.15.03.12}. 
We will show that they are equivalent for $n \leq 3$. 
Its $\Q$-part is the Freeness Conjecture from \cite{G1}-\cite{G2} 
under a scissor congruence hat.

The $\varepsilon$--part of Conjecture \ref{1.12.03.4qq} allows to 
express explicitly 
the $F^*$--modules ${\cal E}_n(F)$ via 
the $\Q$--vector spaces ${\cal B}_n(F)$ and $F^*$--modules  $\beta_n(F)$. 

{\bf Examples}. One should have 
$$
{\cal E}_2(F) = \beta_2(F), \quad {\cal E}_3(F) = \beta_3(F)
$$
The $F^*$--module ${\cal E}_4(F)$ should sit in the exact sequence 
$$0 \lra \beta_4(F)\lra  
{\cal E}_4(F) \lra \beta_2(F) \otimes_{\Q} 
{\cal B}_2(F) 
  \lra 0
$$

{\bf 3. Scissor congruence groups and mixed Tate motives over dual numbers.} 
According to 
\cite{Gvol}, Section 1.7,    the Hopf algebra 
$S_{\bullet}(F)$ is isomorphic to the fundamental 
Hopf algebra of the category ${\cal M}_T(F)$ of mixed Tate motives over $F$. 
This means that the category of finite dimensional 
graded comodules over $S_{\bullet}(F)$ is equivalent to the 
category of mixed Tate motives over $F$.


\begin{conjecture} \label{1.15.03.11} 
The category of finite dimensional graded comodules 
over the Lie coalgebra  $
{\cal Q}_{\bullet}(F_{\varepsilon})$ is naturally 
equivalent to a subcategory of the 
category of mixed Tate motives over the dual numbers $F_{\varepsilon}$. 
\end{conjecture}

We show that a  simplex in a Euclidean affine space over $F$ provides 
a comodule over  the Lie coalgebra  $
{\cal Q}_{\bullet}(F_{\varepsilon})$. It corresponds to  
a mixed Tate motive over $F_{\varepsilon}$ obtained by  
perturbation of the zero object over $F$.  
The subcategory ${\cal M}_T(F)$ should be  given by the comodules 
with trivial action of ${\cal E}_{\bullet}(F)$.

 A cycle  approach to the mixed Tate motives 
over $F_{\varepsilon}$ was suggested in 
\cite{BE}.

{\bf The structure of the paper}. 
The additive polylogarithmic motivic complexes 
  are defined in chapter 2. 
In chapter 3 we define  Euclidean scissor congruence 
groups ${\cal E}_n(F)$. In chapter 4 we discuss the category of mixed motives 
over dual numbers and its relationship with the  scissor  congruence groups.

{\bf Acknowledgment}. I am grateful to  Spencer Bloch 
for exiting conversations during my visit to U. 
Chicago at May 2000,  which brought me back to the subject. 
This work was supported by the NSF grant DMS-0099390.

\section{Additive polylogarithmic motivic complexes}

{\bf 1. Cathelineau's complexes}. 
Recall the higher Bloch 
groups ${\cal B}_n(F)$ defined in \cite{G1}-\cite{G2}. One has ${\cal B}_{1}(F) = F^*$.

In \cite{Ca} J-L. Cathelineau defined the $F$--vector spaces, 
denoted
below by $\overline \beta_n(F)$. 
Each of them is generated by the 
elements $<x>_n$, where $x \in F$. One has 
$\overline \beta_1(F) = F$.  
The definition goes by induction. For a set $X$ denote by $F[X]$ the $F$--vector space with the basis $<x>$,  
where $x \in X$. Having the $F$--vector spaces 
$\overline \beta_k(F)$ for $k<n$ one  defines $\overline \beta_{n}(F)$ as the kernel 
of the map of $F$--vector spaces
\begin{equation} \label{1.13.03.591}
\delta: F[F^*-\{1\}] \lra \overline \beta_{n-1}(F) \otimes {\cal B}_{1}(F)
\quad \oplus \quad \overline \beta_1(F) \otimes {\cal B}_{n-1}(F)
\end{equation}
given on the generators by 
\begin{equation} \label{1.12.03.1}
<x> \quad \lms \quad <x>_{n-1} \otimes \{1-x\}_1 \quad + \quad <1-x>_1 \otimes \{x\}_{n-1}
\end{equation}
So, by the very  definition, there is a map of $F$--vector spaces 
\begin{equation} \label{1.13.03.59}
\delta: \overline \beta_{n}(F) \lra \overline \beta_{n-1}(F) \otimes {\cal B}_{1}(F)
\quad \oplus \quad \overline \beta_1(F) \otimes {\cal B}_{n-1}(F)
\end{equation}

Using this let us  define the following complex:
$$
\begin{array}{ccccccccc}
&&\overline \beta_{n-1}(F) \otimes F^* & 
&\overline \beta_{n-2}(F) \otimes \Lambda^2F^* & &&&\\
 \overline \beta_n(F)&\lra &\oplus &\lra &\oplus &\lra     & ... & \lra & F\otimes \Lambda^{n-1}F^*\\
&&F \otimes {\cal B}_{n-1}(F) & &
F \otimes {\cal B}_{n-2}(F) \otimes F^*&&&& 
\end{array}
$$
It has $n$ terms and placed in degrees $[1,n]$. 
Its $k$-th term for $k = 2, ..., n-1$ is 
$$
\overline \beta_{n-k}(F) \otimes \Lambda^kF^* \quad \oplus \quad 
F \otimes {\cal B}_{n-k}(F) \otimes \Lambda^{k-1}F^*
$$
The differential is defined using (\ref{1.12.03.1}) and the Leibniz rule. 
We denote this complex by $\overline {\beta}_{\bullet}(F; n)$.
There is a  homomorphism
$$
F\otimes \Lambda^{n-1}F^* \lra \Omega^{n-1}_{F/\Q}, \qquad a \otimes b_1 \wedge ... \wedge b_{n-1} \lms a \cdot d\log b_1 \wedge ... \wedge d\log b_{n-1}
$$
It provides a homomorphism
$$
H^n \overline {\beta}_{\bullet}(F; n) \lra \Omega^{n-1}_{F/\Q}
$$
It follows from the results of \cite{Ca} that it is an isomorphism. 
Cathelineau conjectured that $\overline {\beta}_{\bullet}(F; n)$ is a 
resolution of 
$\Omega^{n-1}_{F/\Q}[-n]$, i.e. we have 
\begin{conjecture} \label{1.13.03.9}
For $k<n$ one has 
$$
H^k \overline {\beta}_{\bullet}(F; n)\otimes \Q = 0 
$$
\end{conjecture}
It was proved in \cite{Ca} that this is the case for $n=3$.

{\bf 2. Additive polylogarithmic motivic complexes}. 
\begin{definition} \label{1.12.03.7}
The $F^*$--module $\beta_n(F) $  is defined inductively by 
$$
\beta_n(F) := \quad \beta_{n-1}(F)<1> \quad \oplus \quad \overline \beta_n(F) 
$$ 
\end{definition} 
It follows that we have a decomposition
\begin{equation} \label{1.12.03.8}
\beta_n(F) \quad =\quad  \overline \beta_1(F)<n-1> \quad \oplus \quad \overline \beta_{2}(F)<n-2> 
\quad \oplus ... \oplus \quad 
\overline \beta_n(F)
\end{equation}
This  
is the decomposition into eigenspaces of the $F^*$--action.

{\it Example}. The first term in (\ref{1.12.03.8}) is 
$ F<n-1>$. 

We define a complex ${\beta}_{\bullet}(F; n)$ 
just like  $\overline {\beta}_{\bullet}(F; n)$,  but
with $\overline \beta_k(F)$ replaced 
everywhere by $\beta_k(F)$:
$$
\begin{array}{ccccccccc}
&&\beta_{n-1}(F) \otimes F^* & 
&\beta_{n-2}(F) \otimes \Lambda^2F^* & &&&\\
 \beta_n(F)&\lra &\oplus &\lra &\oplus &\lra     & ... & \lra & F\otimes \Lambda^{n-1}F^*\\
&&F \otimes {\cal B}_{n-1}(F) & &
F \otimes {\cal B}_{n-2}(F) \otimes F^*&&&& 
\end{array}
$$

{\it Examples}. 1. The weight two complex ${\beta}_{\bullet}(F; 2)$ 
is 
$$
\beta_{2}(F) \lra F\otimes F^*
$$ 
2. The weight three complex ${\beta}_{\bullet}(F; 3)$ looks as follows:
$$
\begin{array}{ccccc}
&&\beta_{2}(F) \otimes F^* & 
&\\
 \beta_3(F)&\lra &\oplus &\lra & F\otimes \Lambda^{2}F^*\\
&&F \otimes {\cal B}_{2}(F) & &
\end{array}
$$
\begin{proposition} \label{1.12.03.2}
Conjecture \ref{1.13.03.9} for all weights $\leq n$ is equivalent to the following one: 
for $i \leq n$ one has 
\begin{equation}\label{1.12.03.13}
H^{i}{\beta}_{\bullet}(F; n) \quad = \quad \Omega^{i-1}_{F/\Q}<n-i> 
 \end{equation}
\end{proposition}

{\bf Proof}. We have a decomposition into direct sum of complexes
\begin{equation}\label{1.12.03.12}
{\beta}_{\bullet}(F; n)\quad  = \quad \oplus_{k=0}^{n-1}
\overline {\beta}_{\bullet}(F; n-k)<k> 
 \end{equation}
The proposition follows. 

{\bf 3. Additive versus tangential}. 
Recall that the tangent $T{\cal F}$ to a functor ${\cal F}$ from a category of  rings 
to an abelian category is defined by 
$$
T{\cal F}(R):= {\rm Ker}({\cal F}(R[{\varepsilon}]/{\varepsilon}^2) \lra {\cal F}(R))
$$

{\bf Problem}. Show that Suslin's theorem describing the cohomology of 
the Bloch complex remains valid over the dual numbers. 

Using this one could show that the tangent Bloch group 
$TB_2(F)$ is an extension of 
$\beta_2(F)$ by $\Lambda^2F$. 
 
So neither $\overline \beta_n(F)$ nor  $\beta_n(F)$ are 
isomorphic to $T{\cal B}_n(F)$.  
However assuming Conjecture \ref{1.13.03.9} and  thanks to Theorem \ref{1.19.03.10}, 
the complex 
 ${\beta}_{\bullet}(F; n)$ has the same cohomology 
as we expect  for 
 the tangent motivic complex over 
$F$. Moreover it should be quasiisomorphic to it. 
In any case the tangent 
to the polylogarithmic motivic complex (see \cite{G1}-\cite{G2}) should be quasiisomorphic 
to  the complex 
${\beta}_{\bullet}(F; n)$.

 It would be interesting to define the  $F^*$--modules  $\beta_n(F)$ 
via the $F^*$--modules $T{\cal B}_n(F)$, making 
 (\ref{1.12.03.8}) and hence (\ref{1.12.03.12}) theorems.  
A $K$--theoretic definition of 
$\beta_2(F)$ is given by  S. Bloch and H. Esnault in \cite{BE}. 
See also  Sections 3.5-3.6 below for  $n=2,3$.

\section{The Euclidean scissor congruence groups} 

{\bf 1. Euclidean vector spaces}.  
We say that a finite dimensional $F$--vector space has a {\it Euclidean} structure 
if it is equipped with a non-degenerate quadratic form $Q$. 
A {\it Euclidean  affine} space is an affine space over a Euclidean vector space.

A Euclidean structure $Q$ on a vector space $V$ provides a Euclidean 
structure ${\rm det}_Q$ on ${\rm det}V$. A Euclidean volume form in $V$ 
 is a 
volume form ${\rm vol}_Q$   such that 
${\rm vol}_Q^2 = {\rm det}_Q$. Clearly there are two possible choices, 
$\pm {\rm vol}_Q$. A choice of one of them is called an {\it orientation} of $V$. 

Suppose that $V$ is a Euclidean vector space of dimension $2n$. 
Then  a choice of an orientation of  $V$ has the following 
interpretation. 
The Euclidean structure provides an operator 
on  $\ast: \Lambda^{\bullet}V \lra \Lambda^{2n-\bullet}V$ 
such that   $\ast^2=1$. Namely, if 
$ x\in \Lambda^kV$ then for any $y \in \Lambda^{k}V$ one has 
$\ast x \wedge y = <x,y>_Q\cdot {\rm vol}_Q$ where $<>_Q$ is the induced Euclidean 
structure on $\Lambda^kV$. 
The $\ast$--operator leaves invariant the subspace $\Lambda^nV$. 
 Since $\ast^2=1$,  
there is a decomposition 
$$
\Lambda^nV^* = \Lambda^nV^*_+ \oplus \Lambda^nV^*_-
$$ 
on the $\pm 1$ eigenspaces of  $\ast$. 
We call the elements of $\Lambda^nV^*_+$ (respectively  $\Lambda^nV^*_-$) 
the selfdual (respectively antiselfdual) $n$--forms. It follows from the very definition 
that changing the orientation of $V$ we change the $\ast$--operator 
by multiplying it by $-1$, and thus interchange the selfdual and 
antiselfdual $n$--forms.

Let  $F$ be an algebraically closed field. 
Then there  is an alternative geometric description of 
this decomposition. 
 The family of $n$--dimensional isotropic subspaces 
for the  quadratic form  $Q$ has two connected components. 
They are homogeneous spaces for the special orthogonal group of $V$, and interchanged by 
orthogonal transformations with the determinant $-1$. The corresponding isotropic 
subspaces are called the $\alpha$ and $\beta$ planes. 

\begin{lemma} \label{1.18.03.1}
The restriction of any $n$-form from $\Lambda^nV_-^*$ 
(respectively $\Lambda^nV_+^*$)  to 
every isotropic subspace of one (respectively the other) of these families is zero. 
\end{lemma}

{\bf Proof}. Left as an exercise. 

We call the isotropic planes of the first (respectively the second) family the 
$\alpha$-- (respectively $\beta$--planes). Changing the orientation of $V$ 
we interchange the $\alpha$ and $\beta$.

{\bf 2. The Euclidean scissor congruence groups ${\cal E}_n(F)$}. 
We assume that $F$ is an arbitrary field. 
Let $A$ be a Euclidean affine space of dimension $2n-1$. 
A collection of points $x_0, ..., x_{2n-1}$ in $A$ provides 
a simplex with vertices at these points. We say that such a simplex  is 
Euclidean 
if the Euclidean structure in $A$ induces a Euclidean structure 
on each of its faces. 
The abelian group 
${\cal E}_n(F)$ is generated 
by the elements $(x_0, ..., x_{2n-1}; {\rm vol}_Q)$, where $x_i \in A$, 
$(x_0, ..., x_{2n-1})$ is a Euclidean simplex in $A$,  and ${\rm vol}_Q$ is  
a Euclidean volume form. The relations are the following:

i) (Nondegeneracy). $(x_0, ..., x_{2n-1}; {\rm vol}_Q) = 0$ if all $x_i$ belong to a  hyperplane. 

ii) (Skew--symmetry). 
a) $(x_0, ..., x_{2n-1}; -{\rm vol}_Q) = - (x_0, ..., x_{2n-1}; {\rm vol}_Q)$. 

b) For any permutation $\sigma$ one has 
$$
(x_0, ..., x_{2n-1}; {\rm vol}_Q) = {\rm sgn} (\sigma) (x_{\sigma(0)}, ..., 
x_{\sigma({2n-1})}; {\rm vol}_Q) 
$$

iii) (The scissor axiom). For any $n+2$ points $x_0, ..., x_{2n}$ one has 
$$
\sum_{i=0}^{2n}(-1)^i(x_0, ..., \widehat x_i, ..., x_{2n}; {\rm vol}_Q) = 0
$$
provided that all the $2n+1$ simplices involved here are Euclidean. 

iv) (Affine invariance). For any affine transformation $g$ of $A$  one has 
$$
(x_0, ..., x_{2n-1}; {\rm vol}_Q) = 
(gx_0, ..., gx_{2n-1}; g{\rm vol}_{Q}) 
$$
Observe that $g{\rm vol}_{Q}$ is a volume form for the quadratic form $gQ$. 

{\bf Remarks}. 1. This definition makes sense if $2n$ is an integer, 
admitting half integral $n$'s. However the only really interesting Euclidean 
scissor congruence groups 
are the ones in 
odd dimensional spaces, i.e. with $n$ integral. 

2. To get the classical Euclidean scissor congruence groups one has 
to take $F = \R$ and 
 consider only positive definite quadratic forms $Q$ in $\R^{2n-1}$. 
Observe that in this case  
all simplices are Euclidean. 

The dilotations provide an action of the group $F^*$ on 
the group ${\cal E}_n(F)$. 

We define the volume of 
a simplex $S$ spanned by the vectors 
$v_1, ..., v_{n}$ in an $n$--dimensional Euclidean space 
by  $\frac{1}{n!}<v_1\wedge  ... \wedge v_{n}, {\rm vol}_Q>$. 

\begin{lemma} \label{1.18.03.2}
The volume of a simplex provides a homomorphism 
of $F^*$--modules
$$
{\rm Vol}: {\cal E}_n(F) \lra F<n-1>
$$
\end{lemma}

{\bf Proof}. Straitforward. 

{\bf Example}. The length provides an 
isomorphism  ${\cal E}_1(F) \stackrel{\sim}{=} F$. 


{\bf 3. The scissor congruence Hopf algebra $S_{\bullet}(F)$}.  
We assume that $F$ is an arbitrary field. 
The definition given below follows s. 3.4 in  \cite{Gvol}. 
Let $V_{2n}$ is a $2n$--dimensional $F$--vector space and $Q$ 
a non-degenerate quadratic form in $V_{2n}$. Let $M = (M_1, ..., M_{2n})$ 
be a collection of codimension one subspaces in $V_{2n}$. We say 
that $M$ is in generic position to $Q$ if restriction of the quadratic form $Q$ to any 
face $M_I:= M_{i_1} \cap ... \cap M_{i_k}$ is non-degenerate. 

The abelian group $S_n(F)$ is generated by the elements $(M, Q,  {\rm vol}_Q)$, 
usually  denoted simply by $(M, {\rm vol}_Q)$, 
where ${\rm vol}_Q$ is a volume form for the quadratic form $Q$ 
and $M$ is a Euclidean simplex with respect  to $Q$. The relations are the following:

i)  $(M, {\rm vol}_Q) = 0$ if $\cap M_i \not = 0$.

ii) For any  $g \in GL(V_{2n})$ one has $(M, Q, {\rm vol}_Q) = (gM, gQ, g{\rm vol}_Q)$. 

iii) $(M, -{\rm vol}_Q) = -(M, {\rm vol}_Q)$, the 
skew-symmetry with respect to the permutations of $M_i$'s holds.

iv) For any $2n+1$ subspaces $M_0, ..., M_{2n}$ such that for any $I \subset \{0, ..., 2n\}$ 
the restriction 
of $Q$ to $M_I$ is non-degenerate,  and 
$M^{(j)}:= (M_0, ..., \widehat M_j, ..., M_{2n})$,  we get
$$
\sum_{j=0}^{2n} (-1)^j (M^{(j)}, {\rm vol}_Q)=0
$$

{\it Example}. Let $n=1$. Suppose that there is a non zero isotropic vector for $Q$, 
for instance  
$F$ is algebraically closed. Then one has $S_1(F) = F^*$. Indeed, a generator of $S_1(F)$ 
provides an ordered $4$--tuple $(M_1, M_2, L_1, L_2)$ 
of one dimensional subspaces in $V_2$. Here $L_1$ and $L_2$ are the two isotropic subspaces 
for the form $Q$ ordered so that $L_1$ is the $\alpha$--subspace. 
Then the cross--ratio $r(M_1, M_2, L_1, L_2)$  provides an isomorphism of $S_1(F)$ 
with $F^*$.

{\bf Remark}. In \cite{Gvol} we spelled this definition in a bit different form 
by considering algebraic  simplices in the projective space $P(V_{2n})$ 
equipped with a non-degenerate quadric $Q$, and defining an orientation 
by choosing one of the families of maximally isotropic subspaces on the quadric $Q$.
However these two definitions are equivalent. 
Indeed, a choice of the Euclidean volume form ${\rm vol}_Q$ determines the 
$\ast$--operator, and hence,  by Lemma \ref{1.18.03.1}, a 
choice of one of the families of maximally isotropic subspaces on $Q$.  

The commutative, graded Hopf algebra structure on 
$$
S_{\bullet}(F):= \oplus_{n\geq 0} S_{n}(F); \qquad S_{0}(F) = \Q
$$ 
were defined in theorem 3.9 in \cite{Gvol}. 
The space of indecomposables
${\cal Q}_{\bullet}(F)$
of the Hopf algebra $S_{\bullet}(F)$ (see Section 1.2) has a natural structure of a graded Lie coalgebra 
with the cobracket $\delta$ inherited from the coproduct. 

{\bf 4.  The Euclidean Dehn invariant}. It is a homomorphism 
$$
D^E: {\cal E}_{n}(F) \lra \oplus_{k+l=n} {\cal E}_{k}(F)\otimes S_{l}(F), \qquad k,l>0
$$
Let us define its  ${\cal E}_{k}(F)\otimes S_{l}(F)$--component 
 $D^E_{k,l}$.  Choose a partition
$$
\{1, ..., 2n\} = I \cup J; \qquad |I| = 2l
$$
Since $M$ is a Euclidean simplex, $A_I := \cap_{i\in I} M_i$ 
is a Euclidean affine  space of dimension $2k-1$. The hyperplanes $M_j, j\in J$ 
intersect it, providing a collection $\overline M_J$ of $2k$ hyperplanes  there. Choosing 
 a Euclidean 
volume form $\alpha_I$ in $A_I$ we get an element $(\overline M_J,  \alpha_I)\in 
{\cal E}_k(F)$.  

The quotient 
$E_J:= A/A_I$ is a Euclidean vector space of dimension $2l$. 
The hyperplanes $M_i, i\in I$ project to 
the collection of 
hyperplanes  $\overline M_I$ in $E_J$. Choose a volume form $\alpha_J$ 
and let $\alpha =  \alpha_I \otimes \alpha_J$. 
We get an element $(\overline M_J,  \alpha_J) \in S_l(F)$. We set
$$
D^E_{k,l}(M,  \alpha):= \sum_{I} (\overline M_J,  \alpha_I) \otimes (\overline M_I,  \alpha_J)
$$

The statements that $D^E$ is a group homomorphism 
is checked just as Theorem 3.9a)   in \cite{Gvol}. It
follows easily from the very definitions that  
$$
(D^E \otimes {\rm Id} + {\rm Id} \otimes D) \circ D^E =0
$$

Projecting the second component of $D^E$ to ${\cal Q}_l(F)$ we get 
 the reduced Euclidean Dehn invariant 
$$
\overline D^E: {\cal E}_{n}(F) \lra \oplus_{k+l=n} {\cal E}_{k}(F)\otimes {\cal Q}_{l}(F), \qquad k,l>0
$$

\begin{lemma} \label{1.03.15.1} 
a) The Dehn invariant provides ${\cal E}_{\bullet}(F)$ with a structure 
of a comodule over the graded Hopf algebra ${S}_{\bullet}(F)$.

b) The reduced Dehn invariant provides ${\cal E}_{\bullet}(F)$ with a structure 
of a graded comodule over the Lie coalgebra ${\cal Q}_{\bullet}(F)$.
\end{lemma}

{\bf Proof}. a) The proof is similar to the one of Theorem 3.9b) in [G{\rm vol}]. 

b) This is a standard consequence of a). The lemma is proved.

Therefore we get a Lie coalgebra in the category 
$\Q_{\varepsilon}$--mod: 
$$
{\cal Q}_{\bullet}(F_{\varepsilon}):= {\cal Q}_{\bullet}(F) 
\oplus {\cal Q}_{\bullet}(F_{\varepsilon})\cdot {\varepsilon}
$$

The standard cochain complex of the Lie coalgebra 
${\cal Q}_{\bullet}(F_{\varepsilon})$ is given by the complex 
\begin{equation} \label{1.03.15.5}
{\cal Q}_{\bullet}(F_{\varepsilon}) \lra 
\Lambda^2{\cal Q}_{\bullet}(F_{\varepsilon}) \lra ... \lra 
\Lambda^n{\cal Q}_{\bullet}(F_{\varepsilon}) \lra ...
\end{equation}
in $\Q_{\varepsilon}$--mod, where the first map is the cobracket, and the others are defined via the Leibniz rule. 
The decomposition into $\Q$- and ${\varepsilon}$- components 
provides a decomposition 
of this complex into a direct sum of two subcomplexes, called the $\Q$- and ${\varepsilon}$--components. 
It follows from the very definitions that 
the $\Q$--  (respectively ${\varepsilon}$--) component of (\ref{1.03.15.5}) is the 
reduced non-Euclidean (respectively Euclidean) Dehn complex of $F$:
$$
\Lambda^*\Bigl({\cal Q}_{\bullet}(F_{\varepsilon})\Bigr) = 
{\cal Q}^{\ast}(F) \oplus {\cal E}^{\ast}(F)\cdot \varepsilon
$$

\begin{conjecture} \label{1.03.15.3} Suppose that $F$ is an algebraically closed field. 
Then there is canonical isomorphism 
$$
{\Bbb H}^i_{(n)}({\cal Q}_{\bullet}(F_{\varepsilon})) \quad = \quad 
{\rm gr}^{\gamma}_nK_{2n-i}(F)\otimes \Q  \quad \oplus \quad \Omega^{n-i}_{F/\Q}
$$
\end{conjecture}

The $\Q$-part of this isomorphism was conjectured in Section 1.7 in  
\cite{Gvol}. 
The new ingredient is the ${\varepsilon}$-part, 
which is equivalent to Conjecture \ref{1.15.03.12}. 

{\it We will assume from now on that $F$ is an algebraically closed field}.

{\bf 5. The weight two Euclidean Dehn complex}. This is the complex
$$
D^E: {\cal E}_{2}(F) \lra {\cal E}_{1}(F)\otimes S_{1}(F) = F \otimes F^*
$$
We expect that the results of Cathelineau, Dupont, and  Sah imply 
that this complex is canonically isomorphic to the 
additive dilogarithmic complex ${\beta}_{\bullet}(F; 2)$, i.e. there should exist 
canonical isomorphism
$$
l_2: \beta_2(F) \stackrel{\sim}{\lra} {\cal E}_{2}(F)
$$
which commutes with the coproducts and the volume homomorphisms. In other words 
it  makes the following diagram commute:
$$
\begin{array}{ccccc}
F&\stackrel{{\rm vol}}{\longleftarrow} &{\cal E}_{2}(F) &\stackrel{D^E}{\lra} & F \otimes F^*\\
=\downarrow &&\downarrow l_2&&\downarrow =\\
F&\stackrel{{\rm v}_2}{\longleftarrow} &\beta_2(F) &\stackrel{\delta}{\lra} &F \otimes F^*
\end{array}
$$

{\bf 6. The weight three reduced Euclidean Dehn complex}. 
This is the complex
\begin{equation} \label{1.03.15.10}
\begin{array}{ccccc}
&&{\cal E}_{2}(F)\otimes S_{1}(F) &&\\
{\cal E}_{3}(F) & \lra  & \oplus &\lra & F \otimes \Lambda^2F^*\\
&&{\cal E}_{1}(F)\otimes {\cal Q}_{2}(F) &&
\end{array}
\end{equation}
We conjecture that the complex (\ref{1.03.15.10}) is canonically isomorphic to the 
additive trilogarithmic complex ${\beta}_{\bullet}(F; 3)$. 
This means the following. One should have 
canonical isomorphism
$$
l_3: \beta_3(F) \stackrel{\sim}{\lra} {\cal E}_{3}(F)
$$
It should commute with the coproduct and the volume homomorphisms.
Finally, combined with the isomorphism $l_2$, it should   induce an 
isomorphism of complexes
$$
\begin{array}{ccccc}
&&{\cal E}_{2}(F)\otimes S_{1}(F) &&\\
{\cal E}_{3}(F) & \lra  & \oplus &\lra & F \otimes \Lambda^2F^*\\
&&{\cal E}_{1}(F)\otimes {\cal Q}_{2}(F) &&\\
&&&&\\
l_3\downarrow &&\downarrow &&\downarrow =\\
&&&&\\
&& \beta_2(F)\otimes F^* &&\\
\beta_3(F)&\lra  &\oplus&\lra &F \otimes \Lambda^2F^*\\
&&\beta_1(F)\otimes {\cal B}_2(F) &&
\end{array}
$$

{\bf 7. The higher reduced Euclidean Dehn complexes}. 
\begin{conjecture} \label{1.03.15.12} There exist canonical injective homomorphisms 
of $F^*$--modules 
$$
l_n: \beta_n(F) \hra {\cal E}_{n}(F)
$$
which commutes with the coproduct and the volume homomorphisms. 
\end{conjecture}
It follows that the homomorphisms $l_k$ for $k \leq n$ 
provide morphisms of complexes
$$
{\beta}_{\bullet}(F; n) \lra {\cal E}^*(F; n)
$$
One can not expect the maps $l_n$ to be isomorphisms for $n\geq 4$.

\section{The structure of motivic Lie algebras over dual numbers}

{\bf 1. The Tannakian formalism for mixed Tate motives over dual numbers}. 
We expect the category ${\cal M}_T(F_{\varepsilon})$
of mixed Tate motives over $F_{\varepsilon}$ to be a mixed Tate $\Q$--category with 
the Ext groups given by the formula 
\begin{equation} \label{1.20.03.1}
{\rm Ext}^i_{{\cal M}_T(F_{\varepsilon})}(\Q(0), \Q(n)) = {\rm gr}^{\gamma}_n
K_{2n-i}(F_{\varepsilon})_{\Q}, \qquad A_{\Q}:= A \otimes \Q
\end{equation}
The following result 
(Theorem 6.5 in \cite{Gvol}) is  deduced from 
Goodwillie's theorem \cite{Go}: 
\begin{theorem} \label{1.19.03.10} For an arbitrary field $F$ one has 
$$
{\rm gr}^n_{\gamma}K_{2n-i}(F_{\varepsilon})\otimes \Q  \quad= \quad 
{\rm gr}^n_{\gamma}K_{2n-i}(F)\otimes \Q  \quad \oplus  \quad \Omega^{n-i}_{F/\Q}
$$
\end{theorem}

The projection $F_{\varepsilon} \lra F$ and inclusion $F \hra F_{\varepsilon}$ 
give rise to the functors  
$$
{\cal M}_T(F_{\varepsilon}) \lra {\cal M}_T(F); \quad 
{\cal M}_T(F_{\varepsilon}) \lra {\cal M}_T(F)
$$

Therefore the Tannakian formalism implies that there should exist a graded 
Lie coalgebra ${\cal L}_{\bullet}(F_{\varepsilon})$ 
  such that the category ${\cal M}_T(F_{\varepsilon})$
of mixed Tate motives over $F_{\varepsilon}$ is canonically equivalent to the category 
of finite dimensional graded comodules 
over ${\cal L}_{\bullet}(F_{\varepsilon})$. 
Let us denote by $L_{\bullet}(F_{\varepsilon})$ the corresponding Lie algebra. 
It is a semidirect product of its ideal $L^a_{\bullet}(F_{\varepsilon})$ and 
the fundamental Lie algebra $L_{\bullet}(F)$ 
of the category ${\cal M}_T(F)$ of mixed Tate motives over $F$. 
The ideal ${L}^a_{\bullet}(F)$ is not
abelian. 
The group $F^*$ acts by the automorphisms of $F_{\varepsilon}$: 
$\lambda: a+b \varepsilon \lms a+ \lambda b \varepsilon$. 
So it acts by functoriality 
on the  fundamental Lie algebra ${L}_{\bullet}(F_{\varepsilon})$.

The ${\rm Ext}$'s (\ref{1.20.03.1}) can be computed by the weight $n$ part of 
the standard cochain complex of the Lie coalgebra ${\cal L}_{\bullet}(F_{\varepsilon})$. 
So using theorem \ref{1.19.03.10} we get a conjectural formula
\begin{equation} \label{1.13.03.50}
H^i_{(n)}({\cal L}_{\bullet}(F_{\varepsilon})) = 
{\rm gr}^n_{\gamma}K_{2n-i}(F_{\varepsilon})_{\Q}
\end{equation}

Formula (\ref{1.13.03.50}) contains a lot of 
information about the fundamental Lie coalgebra ${\cal L}_{\bullet}(F_{\varepsilon})$. 
For example, it dictates an isomorphism 
\begin{equation} \label{1.13.03.52}
{\cal L}_{1}(F_{\varepsilon}) = F^*_{\varepsilon}\otimes \Q \stackrel{\sim}{=}
F^*_{\Q} \oplus F 
\end{equation}
Indeed $
{\cal L}_{1}(F_{\varepsilon}) = {\rm Ext}^1_{{\cal M}_T(F_{\varepsilon})}(\Q(0), \Q(1)) = K_1(F_{\varepsilon})_{\Q} = F^*_{\varepsilon}\otimes \Q
$. 
Arguing in a similar way we conclude that one 
should have isomorphisms
$$
{\cal L}^a_{2}(F_{\varepsilon}) \stackrel{}{=} T{\cal B}_2(F); \quad {\cal L}^a_{3}(F_{\varepsilon}) \stackrel{}{=} T{\cal B}_3(F)
$$
One should have the canonical injective maps of $F^*$--modules
$$
T{\cal B}_n(F) \hookrightarrow {\cal L}^a_{n}(F_{\varepsilon})
$$
but they no longer isomorphisms for $n\geq 4$, 
just like in the usual case, see \cite{G2}.

\begin{conjecture} \label{1.12.03.4s} a) There exists canonical 
inclusion of Lie coalgebras
$$
{\cal Q}_{\bullet}(F_{\varepsilon}) \hra {\cal L}_{\bullet}(F_{\varepsilon})
$$

b) It induces an isomorphism 
${\Bbb H}^i_{(n)}({\cal Q}_{\bullet}(F_{\varepsilon})) 
\stackrel{\sim}{\lra} {H}^i_{(n)}({\cal L}_{\bullet}(F_{\varepsilon}))$.

c) The Lie subcoalgebra ${\cal Q}_{\bullet}(F_{\varepsilon})$ 
is characterized by a) and b).
\end{conjecture}

The part a) is nothing else but reformulation of Conjecture \ref{1.15.03.11}. 
Theorem \ref{1.19.03.10} shows that the $\varepsilon$-part of b) 
is equivalent to Conjecture \ref{1.15.03.12}.

{\bf 2. The strong version of the Freeness Fonjecture \ref{1.12.03.4qq}}.  
Since 
$$
\frac{Q_{\bullet}(F_{\varepsilon})}{I_{\bullet}(F_{\varepsilon}) } = 
Q_{-1}(F_{\varepsilon}), \qquad H_1(I_{\bullet}(F_{\varepsilon})) = 
\frac{I_{\bullet}(F_{\varepsilon})}{[I_{\bullet}(F_{\varepsilon}), I_{\bullet}(F_{\varepsilon})]}
$$ 
we get an action 
$
Q_{-1}(F_{\varepsilon})\otimes  H_1(Q_{\bullet}(F_{\varepsilon}))
\lra H_1(Q_{\bullet}(F_{\varepsilon}))
$. 
Dualizing it and using (\ref{1.13.03.52}) we come to the map
\begin{equation} \label{1.13.03.57}
H^1(I_{\bullet}(F_{\varepsilon}))\lra H^1(I_{\bullet}(F_{\varepsilon})) \otimes \Bigl(F^*_{\Q} \oplus F \cdot \varepsilon\Bigr)
\end{equation} 
According to Conjecture \ref{1.12.03.4qq} the degree $-n$ part 
of the 
$\varepsilon$--component of the map (\ref{1.13.03.57}) can be identified 
with the map
$$
\beta_n(F) \quad \lra \quad 
{\cal B}_{n-1}(F) \otimes F^*_{\Q} \quad \oplus \quad \beta_{n-1}(F) \otimes F \cdot 
\varepsilon 
$$
We strengthen  Conjecture \ref{1.12.03.4qq} by adding to it  
that this map coincides with  (\ref{1.13.03.59}) after interchanging the 
 factors in the second term.

\begin{proposition} \label{1.12.03.14}
The strong version of conjecture \ref{1.12.03.4qq} is equivalent to formula (\ref{1.12.03.13}) 
for all $n$, and hence to conjecture 
\ref{1.13.03.9}.
\end{proposition}

{\bf Proof}. The same argument using the Hochshild-Serre spectral sequence for the ideal 
$I_{\bullet}(F_{\varepsilon})$ as in \cite{G2} works.  The second statement follows from 
Proposition \ref{1.12.03.2}.

Address: Dept. of Math. Brown University, Providence RI 02912, USA.

e-mail: sasha@math.brown.edu
\end{document}